# Thermal modelling for endocardiac radiofrequency ablation: comparison of hyperbolic bioheat equation and Pennes bioheat equation with finite element method

**Pengfei Liu**[1,2]**, Jiquan Liu**[1,2] **and Huilong Duan**[1,2]
[1]College of Biomedical Engineering & Instrument Science, Zhejiang University
[2]Key Laboratory for Biomedical Engineering, Ministry of Education, China
Hangzhou, Zhejiang Province, P.R. China, 310027

E-mail: liujq@zju.edu.cn

**Abstract.** The objectives of this study are to model the endocardiac radiofrequency (RF) ablation procedure and to employ the Hyperbolic Bioheat Equation (HBE), which takes the thermal wave behaviour into account, comparing the results with those obtained using the common Pennes Bioheat Equation (BE) method. A complex model is created to cover particular endocardiac physical and geometry environment. Finite Element Method (FEM) is adopted to study the model with both BE and HBE methods. Different convection coefficients and voltages are applied to simulate different conditions. Lesion size, max temperature and specified position temperature are selected as criteria to evaluate the simulated results. The study found that during ablation, the lesion size difference ratio can reach 20% in some periods. The difference is obvious and cannot be neglected.

Keywords- Finite Element Method, hyperbolic heat equation, Pennes heat equation, , radiofrequency ablation, theoretical model

## 1. Introduction

Radiofrequency (RF) heating of biological tissues is currently applied in many surgical and therapeutic procedures such as the elimination of cardiac arrhythmias, the destruction of tumours, the treatment of gastroesophageal reflux disease, and the heating of the cornea for refractive surgery(Berjano, 2006). In order to investigate and develop new RF ablation techniques, apart from understanding the complex electrical and thermal phenomena involved in the heating process, many theoretical models have been employed(Schutt and Haemmerich, 2008; Lai *et al.*, 2004; Berjano and Hornero, 2004; Tungjitkusolmun *et al.*, 2000; Tungjitkusolmun *et al.*, 1997). Currently, in most of the models, Pennes Bioheat Equation (BE) is applied(Pennes, 1948), in which the heat conduction term is based on Fourier's theory. It relates to heat flux ($\vec{q}$) in the following way:

$$\vec{q}(x,t) = -k\nabla T(x,t) \tag{1}$$

where T(x, t) is the temperature at point x in the domain at the time t, and k is the thermal conductivity. This theory assumes an infinite thermal energy propagation speed, i.e. any local temperature disturbance causes an instantaneous perturbation in the temperature at each point in the medium (Jing et al., 1999). For most heating processes, this assumption is suitable.

In some situations, however, such as very low temperature, very high heat flux or very short heating duration, the Fourier's heat conduction theory breaks down. In those conditions, the wave nature of heating processes is more pronounced(Shih et al., 2005). Vick and Özisik investigated the wave characteristics of heat propagation in a semi-infinite medium containing volumetric energy sources(Vick and Ozisik, 1983), and found that thermal behaviour in some situations cannot be predicted by the classical linear or nonlinear diffusion theory because it allows for the immediate diffusion of heat as soon as the energy is released, without considering the effect due to a relaxation time or start-up time. In addition, Mitra found that the wave nature of heat transfer is slow enough in processed meat(Mitra et al., 1995). To take this phenomenon into account, both Cattaneo and Vernotte suggested a modified heat flux model in the following form:

$$\vec{q}(x, t+\tau) = -k\nabla T(x,t) \qquad (2)$$

where $\tau$ is the thermal relaxation time. This equation assumes that the effect (heat flux) and the cause (temperature gradient) occur at different times and the delay between heat flux and temperature gradient is $\tau$. That is, the heat flow does not start instantaneously, but grows gradually with a thermal relaxation time $\tau$, after applying a temperature gradient(Shih et al., 2005).

This equation is normally called Hyperbolic Bioheat Equation (HBE). In the particular case where $\tau=0$, it is obviously equivalent to the Fourier theory.

The wave nature of heat transfer in living tissue may play an important role during rapid heating, such as during thermal ablation/thermal surgery when using high-intensity focused ultrasound or radiofrequency ablation. HBE method is commonly adopted in very short heating durations such as during laser heating. For radiofrequency ablation, the equation is rarely applied. In this study, a complex model is created to cover particular endocardiac physical and geometry environments. The comparative study investigates the BE and HBE methods to find the differences.

## 2. Model and Method

*2.1. Analytical approach*

Considering the finite thermal propagation speed, Cattaneo and Vernotte formulated a modified unsteady heat conduction equation as follows:

$$-\nabla q(x,t) + Q(x,t) = \rho c \frac{\partial T(x,t)}{\partial t} \qquad (3)$$

where $\rho$ is the density (kg/m$^3$) and $c$ is the specific heat (J/kg·K), with the heat transfer model derived from Equation (2) in more detailed format:

$$\vec{q}(x,t) + \tau \frac{\partial \vec{q}(x,t)}{\partial t} = -k\nabla T(x,t) \qquad (4)$$

We can get the following equation:

$$\frac{\partial^2 T}{\partial t^2}(x,t) + \frac{1}{\tau}\frac{\partial T}{\partial t}(x,t) - \frac{k}{\rho c \tau}\Delta T(x,t) =$$
$$\frac{1}{\rho c \tau}Q(x,t) + \frac{1}{\rho c}\frac{\partial Q}{\partial t} \tag{5}$$

Shih(Shih et al., 2005) found that the thermal relaxation time of tissues will cause delay of the appearance of the peak temperature during thermal treatments. The lag behaviour of the peak temperature would result in a lower thermal dose level. For rapid heating (i.e. the heating duration is shorter than the thermal relaxation time of tissue τ), the temperature predicted by the BE is higher than that of the HBE. For a rapid heating process, the HBE may provide an appropriate way of describing the dimensions of thermal lesion during thermal treatments. In contrast, the difference between the dimensions of thermal lesion predicted by the BE and the HBE will become smaller when the heating duration is longer than the thermal relaxation time of tissue.

Taking into consideration the real endocardiac RF ablation environment, the blood convection cooling effect plays an important role when the lesion result is not clear.

*2.2. Specified endocardiac RF ablation model*
The HBE method is commonly adopted in very short heating duration cases such as in laser heating. For radiofrequency ablation, the equation is rarely applied. Berjano (Romero-Garcia *et al.*, 2009; Trujillo M, 2009; Molina *et al.*, 2010; Rivera *et al.*, 2010; Molina *et al.*, 2008a, 2009; Molina *et al.*, 2008b) has studied this topic and proposed a simple model to compare the different results between BE and HBE with both the Finite Element Method and the analytical approach.

Berjano's model used a spherical dry electrode completely embedded in the biological tissue and a constant-power source was adopted to simplify the electrical thermal effect. This model is not suitable for the complex endocardiac radiofrequency ablation environment for the following reasons:

- Berjano's model only covers blood perfusion effect by capillary vessel, and the blood flow convection in the heart chamber is not involved.

- In Berjano's model, biological tissue has an infinite dimension but the real heart chamber wall's thickness is less than 2cm.

- In Berjano's model, a constant power source is applied. But in an actual ablation generator, temperature control and power control mathematics are used. They are both fulfilled with voltage adjustment.

Considering the particular endocardiac environment, the authors created another model, shown in Figure 1. The model is similar to Tungjitkusolmun's work(Tungjitkusolmun et al., 1997), but its shape is a cylinder instead of a cube. The ablation electrode is 5 mm long, and its diameter is 2.6 mm. The thickness of the cardiac tissue is 8 mm and its radius is 20 mm. The cardiac muscle piece lies in the middle of the blood pool and the blood pool depth is 32 mm. The tip of the electrode extends 1.3 mm into the tissue. The depth of the overall model is 40 mm, excluding the electrode.

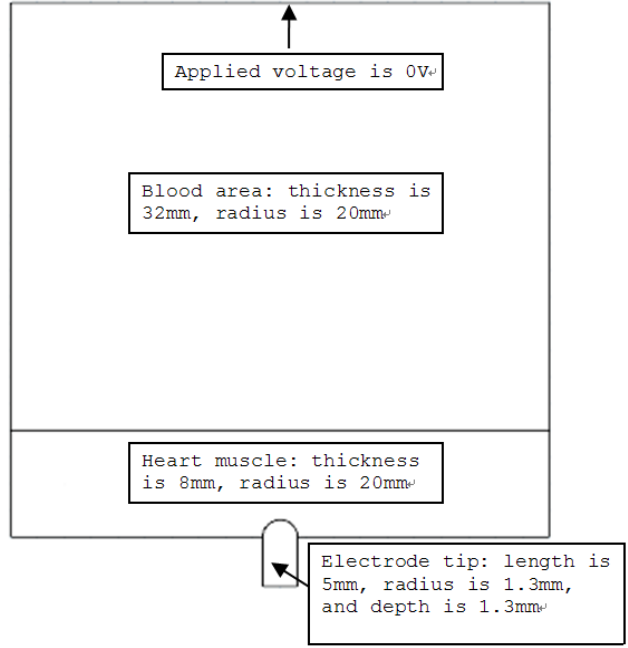

Figure 1.  Schematic diagram of the ablation model geometry

## 2.3. Physical Parameters

In the endocardiac RF ablation model, there are 3 different components: electrode, heart muscle and blood. The 3 different materials' physical properties are listed in Table I. The values for the properties are derived from the fundamental research(Incropera and DeWitt, 1985; Rabbat, 1990) and have been adopted by Tungjitkusolmun(Tungjitkusolmun et al., 1997) .

Table 1.  Physical properties of different regions

| Region | $\rho$ [kg/ m$^3$] | c [J/kg·K] | k[W/m·K] | $\sigma$ [S/m] |
| --- | --- | --- | --- | --- |
| Electrode | 21500 | 132 | 71 | $4\times10^6$ |
| Heart muscle | 1200 | 3200 | 0.550 | 0.222 |
| Blood | 1000 | 4180 | 0.543 | 0.667 |

There is a lack of experimental data on the value of the thermal relaxation time ($\tau$). Values ranging from 10 to 50 s have been found for non-homogeneous inner structure materials(Kaminski, 1990), A value of $\tau$= 16 s has been measured in processed meat(Mitra et al., 1995) which has been adopted in some studies(Lopez Molina et al., 2009; Shih et al., 2005), so in the model, $\tau$ is set to 16 seconds as well.

From the electrical aspect, the Dirichlet boundary conditions are applied. The voltages on the outer surfaces of the model are assumed to be 0 V. The voltage power is applied on the surface of the ablation electrode. The temperature of the blood is assumed to be constant at 37$^o$C due to high flow in the heart chamber.

To simulate the effect of blood convection, a convective boundary condition, he, is assigned at the interface between the electrode and the blood. A convective boundary condition, hc, for the surfaces touching the blood in the chamber is also assigned.

In order to evaluate the difference between the two heat equation methods, various voltages and convection coefficients are adopted to cover the various physical environments. In this simulation, they were divided into 2 different groups.

Group one: Voltage is fixed at 30V. The convection coefficients are: he= 2000 J/( $m^3$·s·K) and hc = 40 kJ/ ($m^3$·s·K) , then various ratios are adopted to simulate various convection coefficients: 0.0, 0.1, 0.25, 0.5, 1.0 and 1.5.

Group two: Convection coefficients are fixed as the standard just like in group one. Different voltages are applied: 25V, 30V, 35V and 40V.

Some instances, such as convection coefficient ratio of 0, do not happen in the real environment, but looking at the trends among different parameter values can provide useful detailed information.

*2.4. Use of COMSOL Multiphysics*

In the model, there are two different physical procedures: electric heating and heat transfer. Electric heating exists in both the heart muscle and blood area. For heat transfer, there are a total of 5 components:

- heat transfer in the heart muscle

- heat transfer in the blood area

- heat transfer in the electrical tip

- heat convection cooling on the tip surface

- heat convection cooling on the heart muscle surface

Both BE and HBE methods should only be considered for heat transfer in the heart muscle. For the other four parts, the standard heat procedure is applied.

In this study, the FEM simulation was performed using COMSOL Multiphysics software version 4.0. COMSOL Multiphysics is an engineering, design, and finite element analysis software environment for the modelling and simulation of any physics-based system. It is a finite element analysis, solver and Simulation software / FEA Software package for various physics and engineering applications, especially coupled phenomena, or multiphysics.

COMSOL is able to fulfil the BE method with the standard heat transfer module, and the method can be accomplished with selection of the "electric current" module and "heat transfer" module. For the HBE method, there is no existing module in COMSOL. But COMSOL can fulfil any equation by user-defined PDE if its equation format is clear, so the "electric current" module and "PDE" module were selected. The "coefficient form PDE" option was used to do a time-dependent analysis. The time interval step was set to 0.1 second and end time was 120 seconds.

In order to verify this method, both 2D axial symmetry model and 3D model were created. In each model, 3 different heat transfer methods were accomplished:

1) Standard COMSOL BE heat transfer method

2) PDE BE heat transfer method

3) PDE Hyperbolic heat transfer method

The results in method 1 and method 2 were found to be the same, with the 3D model result being equivalent to the 2D result as well. This ensures that the results are reliable. In the paper, the results are all from the 2D axial symmetry model, so only half of the model has been created. Mesh mode is set to "Physics-controlled mesh" and element size is set to "fine". There are a total of 2882 nodes in the 2D model.

## 3. Results and Discussion

FEM analysis was carried out on the model with both BE and HBE methods and the different results were compared with the following criteria: temperature profile, lesion area size, maximum temperature and temperature plot at specified location. Nath et al(Nath *et al.*, 1993) found that when the tissue temperature reaches 50°C, irreversible myocardial injury occurs during RF cardiac ablation, so it is selected as the criteria for lesion size measurement.

*3.1. Voltage is 30V, convection coefficient ratio is 1.0*
Figures 2 and 3 are the temperature profile map at 120 seconds when voltage is 30V and convection coefficient ratio is 1.0. The maximum temperature position lies in the deeper position of the electrode. The red contour line indicates the position where temperature equals 50 °C, and the area covered by the red contour line is the lesion region. The temperature profile shape results are similar between BE and HBE methods.

The lesion area size progress is listed in Figure 4. The solid line is the result using the BE method and the dotted line is the result using the HBE method. The lesion size with the HBE method is less than that of the BE method in the first 30 seconds, following which it surpasses the BE method.

The maximum temperature progress is plotted in figure 5. The maximum temperature in the HBE method is less than that in the BE method in the first 20 seconds, following which it surpasses the BE method.

Figure 6 plots the temperature progress with time in 3 different positions and depths of 1.3mm, 2.6mm and 5.2mm. In this model, the electrode intersection depth is 1.3mm, which indicates the position on the electrode surface. The maximum temperature was found to be in the position of 3.2mm, out of the 3 locations. The temperature with the HBE method is less than the relative one using the BE method in the first period. Following that, HBE surpasses BE in the different times. In the position of 1.3mm, the temperature remains constant throughout most of the ablation time.

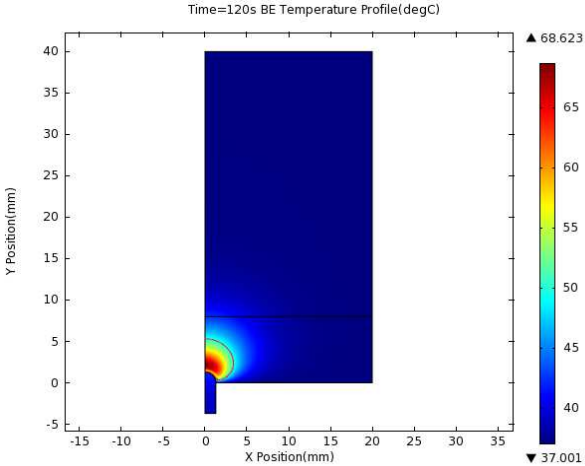

Figure 2. BE temperature profile at 120 seconds

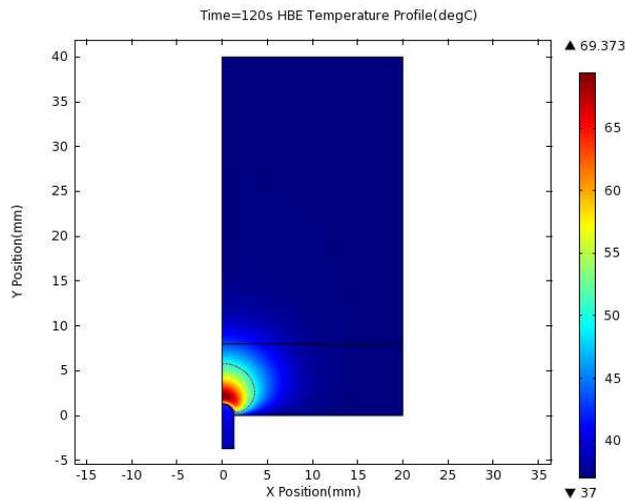

Figure 3.  HBE temperature profile at 60 seconds

During RF ablation heating, the Joule heat is concentrated in the vicinity of the electrode tip. The maximum temperature should exist on the tip surface (1.3mm). Under the influence of the electrode and blood convection cooling, the heat is transferred to the electrode and blood, so its temperature decreases and the max temperature appears in the deeper position of the electrode tip. Comparing between the BE and HBE methods, it can be observed that there is a cross point in lesion size, max temperature value and temperature in the specified location.  As predicted by Shih(Shih et al., 2005), the result is different. It is caused by the blood convection cooling effect.

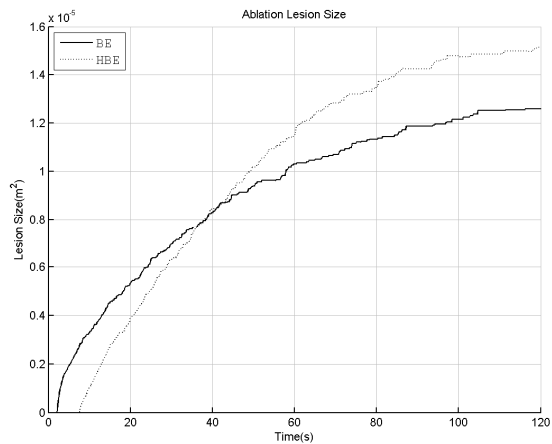

Figure 4.  Lesion size progress plot in 120 seconds

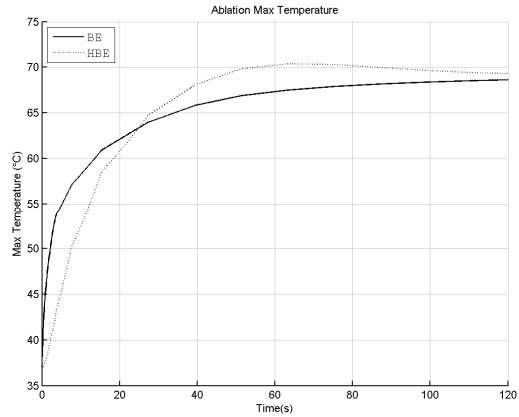

Figure 5.   Max temperature progress plot with time

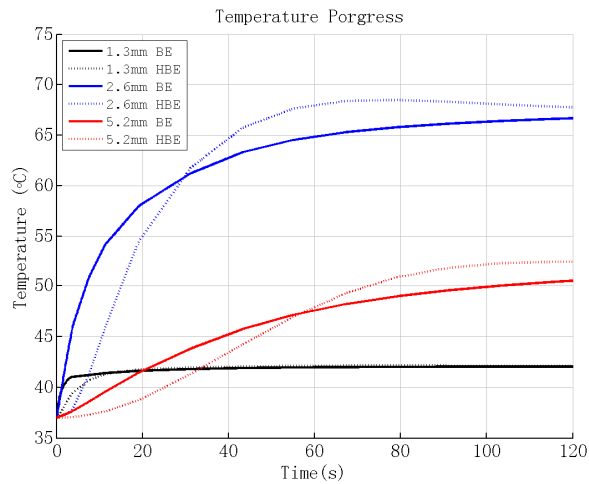

Figure 6.   Temperature progress in the y positions 1.3mm, 2.6mm and 5.2mm.

### 3.2. Various convection coefficient ratios result

Figure 7 shows the lesion size trends in all convection coefficient ratios. The maximum temperature progress is shown in figure 8. The BE method ablation size is always larger than the HBE method when the ratio is 0, but the difference decreases with time. For other cases, their trends are similar: the BE ablation lesion size is larger than the HBE method, after which HBE surpasses BE. The cross point occurs earlier with the increase of the ratio. With increase of the ratio, the lesion size decreases.

For max temperature, the BE is initially larger than the HBE method, then a cross point occurs and the HBE becomes larger than the BE method. The cross point occurs in all ratio cases. When the ratio is 0, the max temperature overruns 100$^{\circ}$C. However, this finding is invalid because the ratio does not reach 0 in the heart.

### 3.3. Various voltages result

Figure 9 shows lesion size trends plotted in different voltages. And the maximum temperature progress is shown in figure 10. All trend lines are similar to the case when the voltage is 30V. The cross point always occurs. With the increase of voltage, both lesion size and max temperature increase.

The cross point always occurs in all cases. With the increase of voltage, the cross point occurs as time decreases.

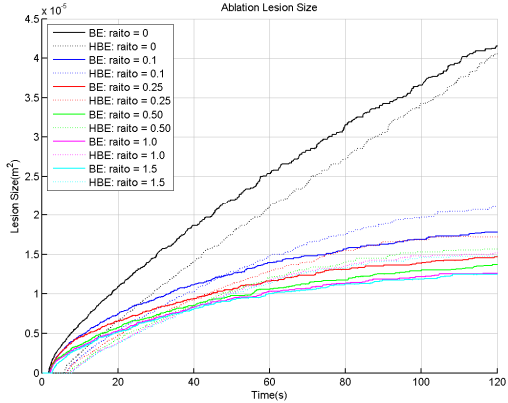

Figure 7. Lesion size progress with different convection ratios

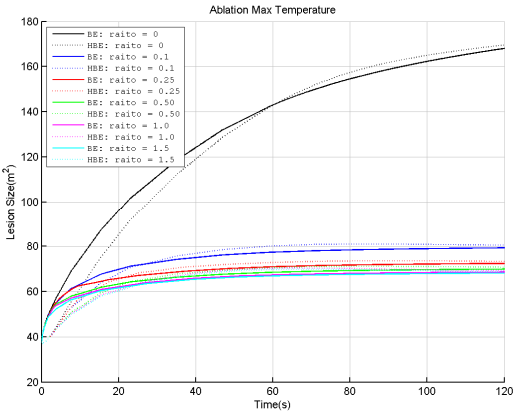

Figure 8. Max temperature progress with different convection ratios

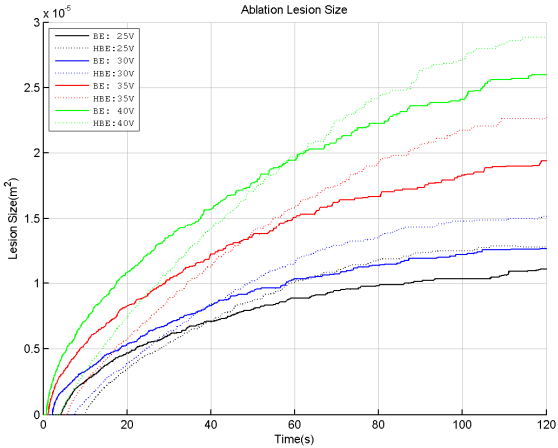

Figure 9. Lesion size progress with different voltages

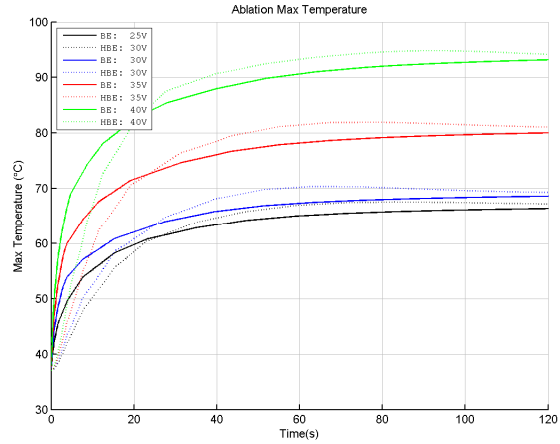

Figure 10. Max temperature progress with different voltages

*3.4. Analysis*

Figures 11 and 12 show the lesion size difference ratio between BE and HBE method in various convection ratios and voltage cases. In the first few seconds, the lesion size is small, and it is useless to evaluate the lesion size difference ratio. After 30 seconds, the lesion size difference ratio varies with time. When the convection coefficient ratio is not 0, the lesion size difference ratio can reach 20% in some periods. The difference is obvious and cannot be neglected.

Figure 13 shows the heart energy absorption progress in different convection coefficient ratios. When the ratio is 0, most of the power is absorbed by the heart tissue and the BE method gets more energy than the HBE method. The lesion size result is the same, as predicted by Shih(Shih et al., 2005). Heart tissue receives more heat energy with BE method that that with HBE method, so its lesion size is always larger. When the ratio is not 0, the heart muscle does not absorb most of the energy, and blood flow takes over more than half of the energy. With the increase of ratios, more energy is transferred into the blood flow. As the ablation time increases, the heart initially absorbs more heat energy with the BE method, after which the HBE method surpasses the BE method. This is due to the fact that the HBE method limits the heat wave propaganda speed, so it can prevent heat from flowing into the blood.

Figure 14 plots the heat absorption by different tissues when the convection coefficient ratio is 1. It is observed that more than half of the energy is absorbed by blood flow, and heart tissue absorption takes the second position.

Therefore, it can be concluded that when RF ablation begins to heat in HBE mode, the flux transfer is limited with specified speed, so its lesion size is smaller than in BE mode. With the increase of time, more energy is transferred into the blood in BE mode. The HBE method produces more energy than the BE method. Therefore, its max temperature and lesion size will surpass that of the BE method. When there is no convection effect, the heart tissue absorbs more energy in the BE method than in the HBE method, so its lesion size is larger than with the HBE method.

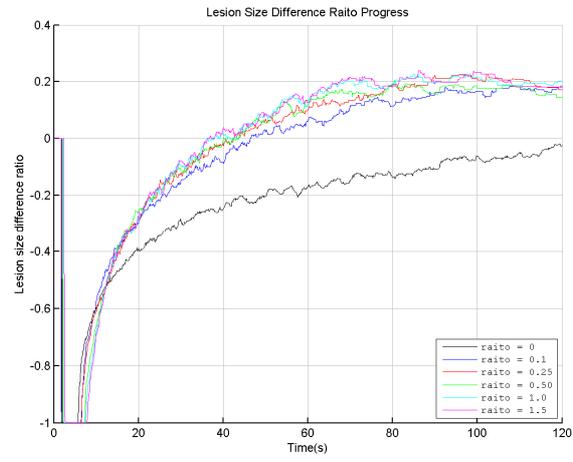

Figure 11. Lesion size difference ratio progress with different ratio

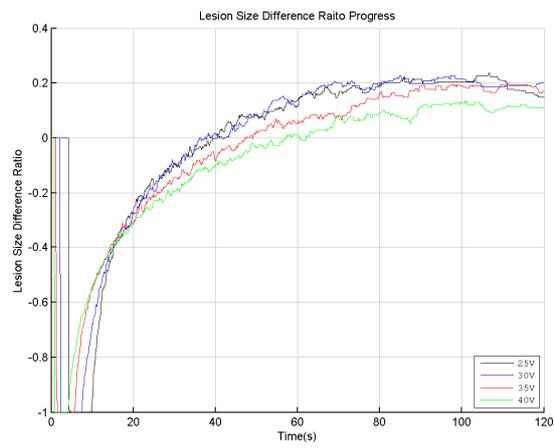

Figure 12. Lesion size difference ratio progress with different voltages

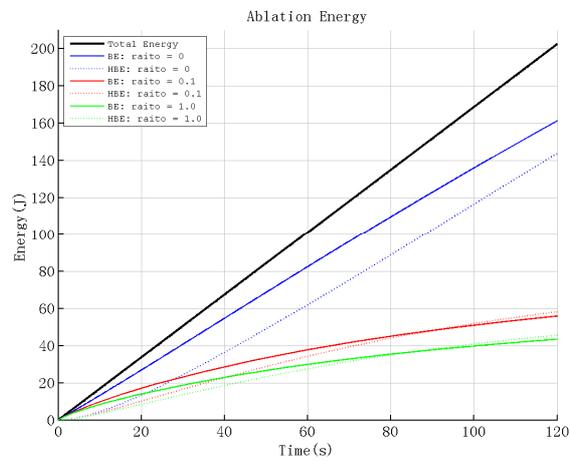

Figure 13. Heart energy absorption progress with different ratios

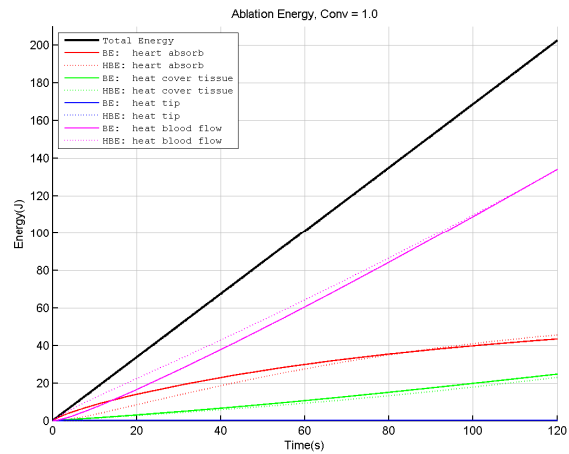

Figure 14. Total energy plot when convection ratio is 1

## 4. Conclusion

In this paper, the authors studied the endocardiac RF ablation procedure using both BE and HBE methods. At first, a specified endocardiac RF ablation model is created, and then the two methods were carried out in various convection coefficient and voltage conditions using FEM software COMSOL. Finally, the different results were compared. It was found that the difference in maximum temperature and lesion area was obvious and cannot be neglected within the specific period of 120 seconds. Considering the more complex properties such as variable material properties with temperature, the variable relaxation time and electric tip category, this comparison should be carried out in a more detailed environment.

**Acknowledgements**

The authors would like to acknowledge National Natural Science Foundation of China for supporting through grant number 30800251 and the Fundamental Research Funds for the Central Universities of China. The authors would also like to thank Prof. Kwan-Hoong Ng for his contributions to this work.

**References**
Berjano E J 2006 Theoretical modeling for radiofrequency ablation: state-of-the-art and challenges for the future *Biomed. Eng. Online* **5**
Berjano E J and Hornero F 2004 Thermal-electrical modeling for epicardial atrial radiofrequency ablation *IEEE Trans Biomed Eng* **51** 1348-57
Incropera F P and DeWitt D P 1985 *Fundamentals of Heat and Mass Transfer* (New York: Wiley)
Jing L, Xu C and Xu L X 1999 New thermal wave aspects on burn evaluation of skin subjected to instantaneous heating *Biomedical Engineering, IEEE Transactions on* **46** 420-8
Kaminski W 1990 Hyperbolic Heat Conduction Equation for Materials With a Nonhomogeneous Inner Structure *Journal of Heat Transfer* **112** 555-60
Lai Y C, Choy Y B, Haemmerich D, Vorperian V R and Webster J G 2004 Lesion size estimator of cardiac radiofrequency ablation at different common locations with different tip temperatures *IEEE Trans Biomed Eng* **51** 1859-64
Lopez Molina J A, Rivera M J, Trujillo M and Berjano E J 2009 Thermal modeling for pulsed radiofrequency ablation: analytical study based on hyperbolic heat conduction *Med Phys* **36** 1112-9
Mitra K, Kumar S, Vedevarz A and Moallemi M K 1995 Experimental Evidence of Hyperbolic Heat Conduction in Processed Meat *Journal of Heat Transfer* **117** 568-73

Molina J A, Rivera M J, Trujillo M and Berjano E J 2008a Effect of the thermal wave in radiofrequency ablation modeling: an analytical study *Phys Med Biol* **53** 1447-62

Molina J A, Rivera M J, Trujillo M and Berjano E J 2009 Thermal modeling for pulsed radiofrequency ablation: analytical study based on hyperbolic heat conduction *Med Phys* **36** 1112-9

Molina J A, Rivera M J, Trujillo M, Burdio F, Lequerica J L, Hornero F and Berjano E J 2008b Assessment of hyperbolic heat transfer equation in theoretical modeling for radiofrequency heating techniques *Open Biomed Eng J* **2** 22-7

Molina J A, Rivera M J, Trujillo M, V. R-G and Berjano E J 2010 Including the effect of the thermal wave in theoretical modeling for radiofrequency ablation. In: *XII Mediterranean Conference on Medical and Biological Engineering and Computing,* ed R Magjarevic: Springer Berlin Heidelberg) pp 521-4

Nath S, Lynch C, Whayne J G and Haines D E 1993 Cellular electrophysiological effects of hyperthermia on isolated guinea pig papillary muscle. Implications for catheter ablation. pp 1826-31

Pennes H 1948 Analysis of tissue and arterial blood temperatures in the resting human forearm *Journal of Applied Physiology* **1** 93-122

Rabbat A 1990 Tissue resistivity *Elecrrical Impedance Tomography*

Rivera M J, Lopez Molina J A, Trujillo M, Romero-Garcia V and Berjano E J 2010 Analytical validation of COMSOL Multiphysics for theoretical models of Radiofrequency ablation including the Hyperbolic Bioheat transfer equation *Conf Proc IEEE Eng Med Biol Soc* 3214-7

Romero-Garcia V, Trujillo M, Rivera M J, Molina J A L and Berjano E J 2009 Hyperbolic Heat Transfer Equation for Radiofrequency Heating: Comparison between Analytical and COMSOL solutions. In: *the Proceedings of the COMSOL Conference 2009,* (Milan

Schutt D J and Haemmerich D 2008 Effects of variation in perfusion rates and of perfusion models in computational models of radio frequency tumor ablation *Med Phys* **35** 3462-70

Shih T C, Kou H S, Liauh C T and Lin W L 2005 The impact of thermal wave characteristics on thermal dose distribution during thermal therapy: a numerical study *Med Phys* **32** 3029-36

Trujillo M R M, López Molina JA, Berjano EJ. 2009 Analytical thermal–optic model for laser heating of biological tissue using the hyperbolic heat transfer equation *Mathematical Medicine and Biology* **26** 187-200

Tungjitkusolmun S, EJ W, H C, JZ T, VR V and JG W 2000 Thermal ± electrical finite element modelling for radio frequency cardiac ablation: Effects of changes in myocardial properties *Medical & Biological Engineering & Computing* **38** 562

Tungjitkusolmun S, Hong C, Jang-Zern T and Webster J G 1997 Using ANSYS for three-dimensional electrical-thermal models for radio-frequency catheter ablation. In: *Proceedings of the 19th Annual International Conference of the IEEE,* pp 161-4

Vick B and Ozisik M N 1983 Growth and Decay of a Thermal Pulse Predicted by the Hyperbolic Heat Conduction Equation *Journal of Heat Transfer* **105** 902-7